\documentclass[10pt,a4paper]{amsart}
\usepackage{amssymb,amsmath}

\headsep=1cm \numberwithin{equation}{section}
\newtheorem{theorem}{Theorem}[section]
\newtheorem{lemma}[theorem]{Lemma}
\newtheorem{proposition}[theorem]{Proposition}
\newtheorem{corollary}[theorem]{Corollary}

\theoremstyle{definition}
\newtheorem{definition}[theorem]{Definition}

\newtheorem{setup and notation}[theorem]{Setup and notation}
\newtheorem{example}[theorem]{Example}

\newcommand{\Gfd}{\operatorname{Gfd}}
\newcommand{\Gid}{\operatorname{Gid}}
\newcommand{\Gpd}{\operatorname{Gpd}}
\newcommand{\G}{\operatorname{G}}
\newcommand{\Ed}{\operatorname{Ed}}
\newcommand{\pd}{\operatorname{pd}}
\newcommand{\fd}{\operatorname{fd}}
\newcommand{\FID}{\operatorname{FID}}

\newcommand{\id}{\operatorname{id}}
\newcommand{\Ext}{\operatorname{Ext}}

\newcommand{\Tor}{\operatorname{Tor}}
\newcommand{\Hom}{\operatorname{Hom}}

\newcommand{\Max}{\operatorname{Max}}

\newcommand{\lo}{\longrightarrow}
\newcommand{\fm}{\frak{m}}

\def\mapdown#1{\Big\downarrow\rlap
{$\vcenter{\hbox{$\scriptstyle#1$}}$}}

\begin{document}

\author[Esmkhani and Tousi]{Mohammad Ali Esmkhani and Massoud Tousi}
\title[Gorenstein homological dimensions and Auslander categories]
{Gorenstein homological dimensions and Auslander categories}

\address{M.A. Esmkhani, Institute for Studies in Theoretical Physics and
Mathematics, P.O. Box 19395-5746, Tehran, Iran-and-Department of
Mathematics, Shahid Beheshti University, Tehran, Iran.}
\email{esmkhani@ipm.ir}

\address{M. Tousi, Institute for Studies in Theoretical Physics and
Mathematics, P.O. Box 19395-5746, Tehran, Iran-and-Department of
Mathematics, Shahid Beheshti University, Tehran, Iran.}
\email{mtousi@ipm.ir}




\begin{abstract} In this paper, we study
Gorenstein injective, projective, and flat modules over a Noetherian
ring $R$. For an $R$-module $M$, we denote by ${\rm Gpd}_RM$ and
${\rm Gfd}_R M$ the Gorenstein projective and flat dimensions of
$M$, respectively. We show that ${\rm Gpd}_RM<\infty$ if and only if
${\rm Gfd}_RM<\infty$ provided the Krull dimension of $R$ is finite.
Moreover, in the case that $R$ is local, we correspond to a
dualizing complex ${\bf D}$ of $\hat{R}$, the classes $A'(R)$ and
$B'(R)$ of $R$-modules. For a module $M$ over a local ring $R$, we
show that $M\in A'(R)$ if and only if ${\rm Gpd}_RM<\infty$ or
equivalently ${\rm Gfd}_RM<\infty$. In dual situation by using the
class $B'(R)$, we provide a characterization of Gorenstein injective
modules.
\end{abstract}

\maketitle

\section{Introduction}

Throughout this paper, $R$ will denote a commutative ring with
nonzero identity and $\hat{R}$ will denote the completion of a
local ring $(R,\fm)$. When discussing the completion of a local
ring $(R,\fm)$, we will mean the $\fm$-adic completion.

Auslander and Bridger [{\bf 3}] introduced the $\G$-dimension,
$\G-\dim_RM$, for every finitely generated $R$-module $M$ (see
also [{\bf 2}]). They proved the inequality $\G-\dim_RM\leq
\pd_RM$, with equality $\G-\dim_RM=\pd_RM$ when $\pd_RM$ is
finite. The $\G$-dimension has strong parallels to the projective
dimension. For instance, over a local Noetherian ring $(R,\fm)$,
the following conditions are equivalent:

(i) $R$ is Gorenstein.

(ii) $\G-\dim_R{R/\fm}<\infty.$

(iii) All finitely generated $R$-modules have finite
$G$-dimension.\\ This characterization of Gorenstein rings is
parallel to Auslander-Buchsbaum-Serre characterization of regular
rings. $\G$-dimension also differs from projective dimension in
that it is defined only for finitely generated modules. Enochs and
Jenda defined in [{\bf 9}] Gorenstein projective modules (i.e.
modules of $\G$-dimension 0) whether the modules are finitely
generated or not. Also, they defined a homological dimension,
namely the Gorenstein projective dimension, $\Gpd_R(-)$, for
arbitrary (non-finitely generated) modules. It is known that for
finitely generated modules, the Gorenstein projective dimension
agrees with the $\G$-dimension. Along the same lines, Gorenstein
flat and Gorenstein injective modules were introduced in [{\bf
9},{\bf 10}].

Let $R$ be a Cohen-Macaulay local ring admitting a dualizing module
$D$. Foxby [{\bf 12}] defined the class $\mathcal{G}_0(R)$ to be
those $R$-modules $M$ such that $\Tor_i^ R(D,M)=\Ext_R^
i(D,D\otimes_RM)=0$ for all $i\geq 1$ and such that the natural map
$M\lo \Hom_R(D,D\otimes_RM)$ is an isomorphism, and
$\mathcal{I}_0(R)$ to be those $R$-modules $N$ such that $\Ext_R^
i(D,N)=\Tor_i^ R(D,\Hom_R(D,N))=0$ for all $i\geq 1$ and such that
the natural map $D\otimes_R\Hom_R(D,N)\lo N$ is an isomorphism. In
[{\bf 11}] Enochs, Jenda and Xu characterize Gorenstein injective,
projective and flat dimensions in terms of $\mathcal{G}_0(R)$ and
$\mathcal{I}_0(R)$.

Let $R$ be a Noetherian ring with dualizing complex ${\bf D}$. The
Auslander categories $A(R)$ and $B(R)$ with respect to ${\bf D}$ are
defined in [{\bf 4}, 3.1]. In [{\bf 5}], it is shown that the
modules in $A(R)$ are precisely those of finite Gorenstein
projective dimension (Gorenstein flat dimension), see [{\bf 5},
Theorem 4.1], and the modules in $B(R)$ are those of finite
Gorenstein injective dimension, see [{\bf 5}, Theorem 4.4]. This may
be viewed as an extension of the results of [{\bf 11}]. Note that,
by [{\bf 4}, Proposition 3.4], if $R$ is a Cohen-Macaulay local ring
with a dualizing module, then an $R$-module $M$ is in $A(R)$ if and
only if $M\in \mathcal{G}_0(R)$ (resp. an $R$-module $M$ is in
$B(R)$ if and only if $M\in \mathcal{I}_0(R)$).

The main aim of this paper is to extend the characterization of
finiteness of Gorenstein dimensions in [{\bf 5}] to arbitrary local
Noetherian rings.

Let $R$ be a local Noetherian ring probably without dualizing
complex, and let ${\bf D}$ denote the dualizing complex of
$\hat{R}$. We define $A'(R)$ to be those $R$-modules $M$ such that
$\hat{R}\otimes_RM\in A(\hat{R})$ and $B'(R)$ to be those
$R$-modules $N$ such that $\Hom_R(\hat{R},N)\in B(\hat{R})$. In
sections 2, 3, and 4, we characterize Gorenstein injective,
projective, and flat modules in terms of the classes $A'(R)$ and
$B'(R)$. To be more precise, we show the following results.

\begin{theorem} Let $R$ be a local Noetherian ring and
$M$ an $R$-module.\\
(i) (See Theorem 2.5) $M$ is Gorenstein flat if and only if $M$
belongs to $A'(R)$ and $\Tor_i^R(L,M)=0$ for all injective
$R$-modules $L$ and all $i>0$.\\
(ii) (See Corollary 3.3) $M$ is Gorenstein projective if and only if
$M$ belongs to $A'(R)$ and $\Ext_R^i(M,P)=0$ for all projective
$R$-modules $P$ and all $i>0$.\\
(iii) (See Theorem 4.8) $M$ is Gorenstein injective if and only if
$M$ belongs to $B'(R)$, $M$ is cotorsion and $\Ext_R^ i(E,M)=0$
for all injective $R$-modules $E$ and all $i>0$.
\end{theorem}

Even more generally, by using the classes $A'(R)$ and $B'(R)$, we
characterize modules of finite Gorenstein injective, projective
and flat dimensions. Namely, we prove the following two results.

\begin{theorem} (See Theorems 3.4 and 3.5) Let $R$ be a Noetherian
ring of finite Krull dimension and $M$ an $R$-module. Then the
following conditions are equivalent:\\
(i) $\Gfd_RM<\infty$.\\ (ii) $\Gpd_RM<\infty$. \\
(More precisely, if $\Gpd_RM<\infty$ or $\Gfd_RM<\infty$, then
$\Max\{\Gfd_RM,\Gpd_RM\}\leq \dim R).$\\
Moreover, if $R$ is local, then the above conditions are equivalent
to the following\\(iii) $M\in A'(R)$.
\end{theorem}

\begin{theorem} (See Theorem 4.10) Let $(R,\fm)$ be a local
Noetherian ring of dimension $d$ and $\Ext_R^i(\hat{R},M)=0$ for
all $i>0$. Then Gorenstein injective dimension of $M$ is finite if
and only if $M$ belongs to $B'(R)$. In particular, if $M\in B'(R)$
then $\Gid_R(M)\leq d$.
\end{theorem}

{\bf Setup and notation} If $M$ is any $R$-module, we use
$\pd_RM$, $\fd_RM$ and $\id_RM$ to denote the usual projective,
flat and injective dimension of $M$, respectively. Furthermore, we
write $\Gpd_RM$, $\Gfd_RM$ and $\Gid_RM$ for the Gorenstein
projective,  Gorenstein flat and Gorenstein injective dimension of
$M$, respectively. Let $\mathcal{X}$ be any class of $R$-modules
and let $M$ be an $R$-module. An $\mathcal{X}$-precover of $M$ is
an $R$-homomorphism $\varphi:X\lo M$, where $X\in \mathcal{X}$ and
such that the sequence,
$$\Hom_R(X',X)\overset{\Hom_R(X',\varphi)}\lo \Hom_R(X',M)\lo 0$$
is exact for every $X'\in \mathcal{X}$. If, moreover, $f
\varphi=\varphi$ for $f\in \Hom_R(X,M)$ implies $f$ is an
automorphism of $M$, then $\varphi$ is called an
$\mathcal{X}$-cover of $M$. Also, an $\mathcal{X}$-preenvelope and
$\mathcal{X}$-envelope of $M$ are defined ``dually''. By $P(R)$,
$F(R)$ and $I(R)$ we denote the classes of all projective, flat
and injective $R$-modules, respectively. Furthermore, we let
$\overline{P(R)}$, $\overline{F(R)}$ and $\overline{I(R)}$ denote
the classes of all $R$-modules with finite projective, flat and
injective dimension, respectively.

We may use the following facts without comment. If $R$ is
Noetherian of finite Krull dimension, then
$\overline{P(R)}=\overline{F(R)}$ (see [{\bf 16}, Theorem 4.2.8
]). Also, if $R$ is Noetherian then for any $M\in
\overline{P(R)}$, we have $\pd_R(M)\leq \dim R$ (see [{\bf 15}, p.
84]).

\section{Gorensein flat dimension}

Let $R$ be a local Noetherian ring an let ${\bf D}$ denote the
dualizing complex of $\hat{R}$. Let $A(\hat{R})$ denote the full
subcategory of ${\bf D}_b(\hat{R})$, consisting of those complexes
$X$ for which ${\bf D}\otimes_{\hat{R}}^{{\bf L}}X\in {\bf
D}_b(\hat{R})$ and the canonical morphism
$$\gamma_X:X\lo {\bf R}\Hom_{\hat{R}}({\bf D},{\bf
D}\otimes_{\hat{R}}^{{\bf L}}X),$$ is an isomorphism. Here, ${\bf
D}_b(\hat{R})$ denote the full subcategory of ${\bf D}(\hat{R})$
(the derived category of $\hat{R}$-modules) consisting of complexes
$X$ with $H_n(X )=0$ for $\mid n \mid
>>0$, see [{\bf 4}].\\
Now, we define $A'(R)$ to be the class of all $R$-modules $M$ such
that $\hat{R}\otimes_RM\in A(\hat{R})$.

\begin{lemma} Let $0\lo M'\lo M\lo M{''}\lo 0$ be an exact sequence of modules
over a local Noetherian ring $R$. Then if any two of $M',M,M{''}$
are in $A'(R)$, so
is the third.\\
\end{lemma}

{\bf Proof.} The exact sequence $0\lo M'\lo M\lo M{''}\lo 0$
yields, the exact sequence $0\lo \hat{R}\otimes_RM'\lo
\hat{R}\otimes_RM\lo \hat{R}\otimes_RM{''}\lo 0$. Now, the
conclusion follows by using [{\bf 5}, Theorem 4.1] and [{\bf 13},
Theorem 2.24]. $\Box$

\begin{proposition}Let $R$ be a local Noetherian ring and let
$M$ be an $R$-module. If $\Gfd_RM<\infty$, then $M\in A'(R)$.
\end{proposition}

{\bf Proof.} By [{\bf 13}, Proposition 3.10], we have
$\Gfd_{\hat{R}}(\hat{R}\otimes_RM)<\infty$. Using [{\bf 5},
Theorem 4.1], we conclude that $\hat{R}\otimes_RM$ belongs to
$A(\hat{R})$. So, the assertion follows by the definition. $\Box$
\\
\\
In the proof of the following lemma we use the method of the proof
of [{\bf 11}, Lemma 3.1].

\begin{lemma} Suppose $K$ is cotorsion of finite flat
dimension and suppose $M$ is an $R$-module. If $\Tor_i^ R(E,M)=0$
for all $i> 0$ and all injective $R$-modules $E$, then
$\Ext_R^i(M,K)=0$ for all $i>0$.
\end{lemma}

{\bf Proof.} We prove by induction on $\fd_RK$. First, let $K$ be
flat and cotorsion. Then $K$ is a summand of a module of the form
$\Hom_R(E,E')$ where $E$ and $E'$ are injective ([{\bf 8}, Lemma
2.3]). It is enough to show that $\Ext_R^ i(M,\Hom_R(E,E'))=0$ for
all $i> 0$. We have $$\Ext_R^i(M,\Hom_R(E,E'))\cong \Hom_R(\Tor_i^
R(M,E),E')$$ for all $i\geq 0$. Thus $\Ext_R^i(M,K)=0$ for all
$i>0$. Now, let $K$ be cotorsion and of finite flat dimension. Let
$F_0\lo K$ be a flat cover of $K$ with kernel $L$. Then $L$ is
cotorsion, see [{\bf 8}, Lemma 2.2]. Also, we have the exact
sequence
$$\Ext_R^i(M,F_0)\lo \Ext_R^i(M,K)\lo \Ext_R^{i+1}(M,L).$$ Since
$K$  and $L$ are cotorsion, then so is $F_0$. Hence, by inductive
hypothesis $\Ext_R^i(M,K)=0$ for all $i> 0$. $\Box$

\begin{lemma} Let $R$ be a  Noetherian ring and $M$ an
$R$-module.\\
(i) If $R$ be a local ring and $M\in A'(R)$, then there
exists a monomorphism $M\lo L$ with $\fd_RL<\infty$.\\
(ii) Assume $\psi:M\lo L$ is a monomorphism such that
$\fd_RL<\infty$ and that $\Tor_i^R(N,M)=0$ for all injective
$R$-modules $N$ and all $i>0$. Then  $M$ possesses a monic
$\overline{F(R)}$-preenvelope $M\lo F$, in which $F$ is flat.\\
(iii) Let $R$-homomorphism $f:M\lo L'$ be an
$\overline{P(R)}$-preenvelope. Assume $\varphi:M\lo L$ is a
monomorphism such that $\pd_RL<\infty$ and that $\Ext_R^ i(M,N)=0$
for all projective $R$-modules $N$ and all $i>0$. Then there exists
a monic $\overline{P(R)}$- preenvelope $M\lo P$, in which $P$ is
projective.
\end{lemma}

{\bf Proof.} (i) Since $M$ belongs to $A'(R)$,
$\Gfd_{\hat{R}}(M\otimes_R\hat{R})$ is finite by the definition
and [{\bf 5}, Theorem 4.1]. Therefore, by [{\bf 5}, lemma 2.19],
we have an exact sequence of $\hat{R}$-modules and
$\hat{R}$-homomorphisms $0\lo M\otimes_R\hat{R}\lo L$, where flat
dimension of $L$ is finite as an $\hat{R}$-module. So, we obtain
an exact sequence $0\lo M\lo L$ of $R$-modules and
$R$-homomorphism, where flat dimension of $L$ is finite as an
$R$-module. Not that every flat $\hat{R}$-module is also flat as
an $R$-module.

(ii) Using [{\bf 7}, Proposition 5.1], there exists a flat
preenvelope $f:M\lo F$. We show that $f$ is
$\overline{F(R)}$-preenvelope. To this end, let $\psi':M\lo L'$ be
an $R$-homomorphism such that $\fd_RL'<\infty$ and let $0\lo K\lo
F'\overset{\pi}\lo L'\lo 0$ be an exact sequence such that
$\pi:F'\lo L'$ is a flat cover. Then $K$ is of finite flat
dimension and also by [{\bf 8}, lemma 2.2], it is cotorsion. Lemma
2.3 implies that $\Ext_R^ i(M,K)=0$ for all $i>0$. So, we have the
exact sequence
$$0\lo \Hom_R(M,K)\lo \Hom_R(M,F')\lo \Hom_R(M,L')\lo \Ext_R^ 1(M,K)=0.$$
Therefore, there exists an $R$-homomorphism $h:M\lo F'$ such that
$\pi h=\psi'$. Since $f:M\lo F$ is flat preenvelope, there exists
an $R$-homomorphism $g:F\lo F'$ such that $h=gf$. Hence, there
exists the $R$-homomorphism $\pi g:F\lo L'$ such that $\pi
gf=\psi'$. Thus $f$ is $\overline{F(R)}$-preenvelope.
Consequently, $f$ is monic, because $\psi$ is monic.

(iii) Since $\varphi:M\lo L$ is monic, it turns out that $f:M\lo L'$
is also monic. Now, let $0\lo K\lo P\overset{\pi}\lo L'\lo 0$ be an
exact sequence such that $P$ is projective $R$-module. It is easy to
see that $K\in \overline{P(R)}$. On the other hand, by hypothesis
and induction on projective dimension, $\Ext_R^ i(M,Q)=0$ for all
$i>0$ and for all $Q\in \overline{P(R)}$. Therefore, $\Ext_R^
i(M,K)=0$ for all $i>0$. Hence $f:M\lo L'$ has a lifting $M\lo P$
which is monic and still an $\overline{P(R)}$-preenvelope. $\Box$

\begin{theorem} Let $(R,\fm)$ be a local Noetherian ring and
$C$ an $R$-module. Then the following conditions are
equivalent:\\
(i) $C$ is Gorenstein flat.\\
(ii) $C$ belongs to $A'(R)$ and $\Tor_i^R(L,C)=0$ for all
injective $R$-modules $L$ and all $i>0$.
\end{theorem}

{\bf Proof.} $(i) \Rightarrow (ii)$ By Proposition 2.2, $C$
belongs to $A'(R)$.
Also, [{\bf 13}, Theorem 3.6], implies the last assertion in (ii).\\
$(ii) \Rightarrow (i)$ By [{\bf 13}, Theorem 3.6], it is enough to
show that $C$ admits a right flat resolution
$${\bf X}=0\lo C\lo F^0\lo
F^1\lo F^2\lo \ldots$$ such that $\Hom_R({\bf X},Y)$ is exact for
all flat $R$-modules $Y$ (i.e. $C$ admits a co-proper right flat
resolution). Lemma 2.4 (i) implies that there exists an exact
sequence $0\lo C\lo L$ of $R$-modules and $R$-homomorphisms such
that $\fd_RL<\infty$. Using Lemma 2.4 (ii), there exists a
monomorphism $f:C\lo K$ which is a flat preenvelope. We obtain the
short exact sequence $0\lo C\overset{f}\lo K\lo B\lo 0$ and so for
every flat $R$-module $F'$ we have the short exact sequence
$$0\lo \Hom_R(B,F')\lo \Hom_R(K,F')\lo \Hom_R(C,F')\lo 0.$$ Let
$E$ be an injective $R$-module. Since $\Hom_R(E,E_R(R/\fm))$ is a
flat $R$-module, we conclude that
$$0\lo C\otimes_RE\lo K\otimes_RE\lo B\otimes_RE\lo 0$$ is an exact
sequence. So, $\Tor_i^ R(E,B)=0$ for all $i>0$ and all injective
$R$-modules $E$, because $K$ is a flat $R$-module. Also, by Lemma
2.1 and Proposition 2.2, we obtain $B\in A'(R)$. Then proceeding
in this manner, we get the desired co-proper right flat resolution
of $C$. $\Box$

\begin{corollary} Let $(R,\fm)$ be a local Noetherian ring of
dimension $d$ and let $M\in A'(R)$. Then
$\Gfd_R(M)=\Gfd_{\hat{R}}(\hat{R}\otimes_RM)$. In particular, if
$M\in A'(R)$ then $\Gfd_RM\leq \dim R$.
\end{corollary}

{\bf Proof.} By [{\bf 13}, Proposition 3.10],
$\Gfd_{\hat{R}}(\hat{R}\otimes_RM)\leq \Gfd_R(M)$. We show that
$\Gfd_R(M)\leq \Gfd_{\hat{R}}(\hat{R}\otimes_RM)$ and so [{\bf
13}, Theorem 3.24] completes the proof. As $M$ belongs to
$A'(R)$, we get that $\hat{R}\otimes_RM$ belongs to
$A({\hat{R}})$. So, by [{\bf 5}, Theorem 4.1]
$\Gfd_{\hat{R}}(\hat{R}\otimes_RM)$ is finite. Set
$\Gfd_{\hat{R}}(\hat{R}\otimes_RM)=t$ and let $$0\lo C\lo
P_{t-1}\lo \ldots\lo P_1\lo P_0\lo M\lo 0,$$ be an exact sequence
of $R$-modules and $R$-homomorphisms such that $P_i^,$s are
projective. We obtain the exact sequence
$$0\lo \hat{R}\otimes_RC\lo \hat{R}\otimes_RP_{t-1}\lo \ldots\lo
\hat{R}\otimes_RP_1\lo \hat{R}\otimes_RP_0\lo \hat{R}\otimes_RM\lo
0.$$ By [{\bf 13}, Theorems 3.14], $\hat{R}\otimes_RC$ is a
Gorenstein flat $\hat{R}$-module. Also, Lemma 2.1 and the above
exact sequence, imply that $C$ belongs to $A'(R)$. In view of
Theorem 2.5, it is enough to show that $\Tor_i^R(C,E)=0$ for all
injective $R$-modules $E$ and all $i> 0$. Let $E$ be an injective
$R$-module and let $\Hom_R(-,E_R(R/\fm))$ denote by $(-)^{\vee}$.
From the natural monomorphism $E\lo (E^{\vee})^{\vee}$, we
conclude that $E$ is a direct summand of $(E^{\vee})^{\vee}$. So,
it is enough to show that $\Tor_i^ R((E^{\vee})^{\vee},C)=0$ for
all $i>0$. By the next result, $\id_{\hat{R}}{(E^{\vee}})^{\vee}$
is finite. It therefore follows from [{\bf 13}, Theorem 3.14] that
$\Tor_i^{\hat{R}}(C\otimes_R{\hat{R}},{(E^{\vee}})^{\vee})=0$ for
all $i>0$. Suppose ${\bf F_{\bullet}}\lo C$ is a flat resolution
of $C$. For every $i>0$, we have
\begin{equation*}
\begin{split}
\Tor_i^ R(C,{(E^{\vee}})^{\vee})& \cong H_i({\bf
F_{\bullet}}\otimes_R{(E^{\vee}})^{\vee})\\&\cong H_i(({\bf
F_{\bullet}}\otimes_R\hat{R})\otimes_{\hat{R}}{(E^{\vee}})^{\vee})
\\&\cong \Tor_i^{\hat{R}}(C\otimes_R{\hat{R}},{(E^{\vee}})^{\vee})
\end{split}
\end{equation*}
The last isomorphism comes from  the fact that
$F_{\bullet}\otimes_R\hat{R}$ is a flat resolution of
$C\otimes_R\hat{R}$, considered as an $\hat{R}$-module. Thus,
$\Tor_i^ R(C,{(E^{\vee}})^{\vee})=0$ for all $i>0$. $\Box$

\begin{lemma} Let $(R,\fm)$ be a local Noetherian  ring
and let $K$ be an $R$-module such that $\id_R(K)$ is finite. Let
$\Hom_R(-,E_R(R/\fm))$ denote by $(-)^{\vee}$. The $R$-module
${(K^{\vee})}^{\vee}$ considered with the $\hat{R}$-module
structure coming from $E_R(R/\fm)$, that is,
$(\hat{r}f)(x)=\hat{r}(f(x))$, for all $\hat{r}\in \hat{R}$, $f\in
\Hom_R(K^{\vee},E_R(R/\fm))$ and $x\in K^{\vee}$. Then
$\id_{\hat{R}}{(K^{\vee})}^{\vee}$ is finite.
\end{lemma}

{\bf Proof.} We deduce that $\fd_R(K^{\vee})$ is finite. It is
easy to see that $\fd_{\hat{R}}(K^{\vee}\otimes_R{\hat{R}})$ is
finite. By the adjoint isomorphism, we have the following
isomorphism
$$\Hom_{\hat{R}}(K^{\vee}\otimes_R{\hat{R}}, E_R(R/{\fm}))\cong
\Hom_R(K^{\vee},E_R(R/{\fm})),$$ as an $\hat{R}$-modules. This
ends the proof, because the injective dimension of
$\Hom_{\hat{R}}(K^{\vee}\otimes_R{\hat{R}}, E_R(R/{\fm}))$ is
finite as an $\hat{R}$-module. $\Box$

\section{Gorensein projective dimension}

In this section, we show that Gorenstein projective dimension of
an $R$-module is finite if and only if its Gorenstein flat
dimension is finite.

\begin{proposition} Let $R$ be a Noetherian ring with finite Krull dimension and
$C$ be an $R$-module. Then $\Gfd_R(C)\leq \Gpd_R(C)$.
\end{proposition}

{\bf Proof.} See [{\bf 13}, Remark 3.3 and Proposition 3.4].
$\Box$

\begin{theorem} Let $R$ be a Noetherian ring of finite Krull dimension
and $M$ an $R$-module. Then the following conditions are equivalent:\\
(i) $M$ is Gorenstein projective.\\
(ii) $\Gfd_RM<\infty$ and $\Ext_R^i(M,P)=0$ for all projective
$R$-modules $P$ and all $i>0$.
\end{theorem}

{\bf Proof.} Assume that $M$ is Gorenstein projective. Then
$\Gfd_RM<\infty$, by Proposition 3.1. Also, [{\bf 13}, Proposition
2.3], implies that $\Ext_R^ i(M,P)=0$ for
all projective $R$-modules $P$ and all $i>0$.\\

Next, we show that $(ii) \Rightarrow (i)$. By [{\bf 13}, Proposition
2.3], it is enough to show that $M$ admits a right projective
resolution $${\bf X}=0\lo M\lo P^ 0\lo P^ 1\lo P^ 2\lo \cdots$$ such
that $\Hom_R({\bf X},Y)$ is exact for every projective $R$-module
$Y$ (i.e. $M$ admits a co-proper right projective resolution).

Since $\Gfd_RM<\infty$, it follows from [{\bf 5}, Lemma 2.19] that
there exists a monomorphism $M\lo L$ with $\fd_RL<\infty$.

Let $i>0$. By assumption and induction on projective dimension,
$\Ext_R^i(M,Q)=0$  for all $Q\in \overline{P(R)}$. On the other
hand, we have
$$\Ext_R^i(M,\Hom_R(E,E'))\cong \Hom_R(\Tor_i^ R(M,E),E')$$for all
injective $R$-modules $E$ and $E'$. Therefore, $\Tor_i^ R(M,E)=0$
for all injective $R$-modules $E$. Note that, for each nonzero
$R$-module $N$, there exists an injective $R$-module $E'$ such that
$\Hom_R(N,E')\neq 0$.

Now, using parts (ii) and (iii) of Lemma 2.4, there exists a
monomorphism $\psi:M\lo Q$ which is a projective preenvelope. We
consider the exact sequence $$0\lo M\overset{\psi}\lo Q\lo B\lo 0,$$
where $B={\rm Coker\; }\psi$. Let $P$ be a projective $R$-module.
Applying the functor $\Hom_R(-,P)$ to the above exact sequence, we
see that $\Ext_R^ i(B,P)=0$ for all $i>0$ because $\psi:M\lo Q$ is a
projective preenvelope.  Also, $\Gfd_RB<\infty$, by [{\bf 13},
Theorem 3.15]. Then proceeding in this manner, we get the desired
co-proper right projective resolution for $M$. $\Box$ \\

We can deduce from Proposition 2.2, Corollary 2.6 and Theorem 3.2
the following result.

\begin{corollary}
Let $R$ be a local Noetherian ring and $M$ an $R$-module. Then
the following conditions are equivalent:\\
(i) $M$ is Gorenstein projective.\\
(ii) $M\in A'(R)$ and $\Ext_R^i(M,P)=0$ for all projective
$R$-modules $P$ and all $i>0$.
\end{corollary}

\begin{theorem} Let $R$ be a Noetherian ring of finite dimension d and
$M$ be an $R$-module. Then the following conditions are
equivalent:\\ (i) $\Gfd_RM<\infty$.\\ (ii) $\Gpd_RM<\infty$.\\
Moreover, if one of the above conditions holds, then $\Gpd_RM\leq
d$.
\end{theorem}

{\bf Proof.} $(i) \Rightarrow (ii)$ We prove the claim by induction
on $\Gfd_RM$. First, let $M$ be a Gorenstein flat $R$-module. Let
$F$ be a flat $R$-module. Consider the minimal pure injective
resolution
$$0\lo F\lo PE^0(F)\lo PE^1(F)\lo \cdots$$ (see [{\bf 16}, pages 39
and 92]). Note that, by [{\bf 16}, Lemma 3.1.6], $PE^n(F)$ is flat
for all $n\geq 0$ and also, by [{\bf 16}, Corollary 4.2.7],
$PE^n(F)=0$ for all $n>d$. Since, every pure injective module is
cotorsion, by [{\bf 13}, Proposition 3.22], $\Ext_R^j(M,PE^i(F))=0$
for all $i\geq 0$ and all $j\geq 1$. Therefore,
$\Ext_R^{d+i}(M,F)\cong \Ext_R^i(M,PE^d(F))$ for all $i\geq 1$, and
so $\Ext_R^{d+i}(M,F)=0$ for all $i\geq 1$. Next, let
$$0\lo C\lo P_{d-1}\lo \cdots\lo P_0\lo M\lo 0$$ be an exact
sequence such that $P_i^,$s are projective. We have
$\Ext_R^{d+i}(M,F)\cong \Ext_R^i(C,F)$ for all $i\geq 1$, and so
$\Ext_R^i(C,F)=0$ for all $i\geq 1$. On the other hand, using [{\bf
13}, Theorem 3.15], we conclude that $\Gfd_RC<\infty$. Therefore, by
Theorem 3.2, $C$ is Gorenstein projective, and hence $\Gpd_RM\leq
d$.

Now, let $\Gfd_RM=t>0$ and let $0\lo K\lo P\lo M\lo 0$ be an exact
sequence such that $P$ is projective. By [{\bf 13}, Proposition
3.12], $\Gfd_RK=t-1$. Hence, induction hypothesis implies that
$\Gpd_RM<\infty$.

$(ii) \Rightarrow (i)$ This follows from Proposition 3.1.

Now, if either ${\rm Gpd}_RM<\infty$ or equivalently
$\Gfd_RM<\infty$, then, by [{\bf 5}, Lemma 2.17], $\Gpd_RM=\pd_RH$,
where $H$ is an $R$-module. This completes the proof. $\Box$ \\

Now, we are ready to deduce the main result of this section by using
Proposition 2.2, Corollary 2.6 and Theorem 3.4.

\begin{theorem} Let $R$ be a local Noetherian ring and $M$
an $R$-module. Then the following conditions are equivalent:\\
(i) $\Gfd_RM<\infty$.\\ (ii) $\Gpd_RM<\infty$.\\(iii) $M\in
A'(R)$.\\ Moreover, if one of the above conditions holds, then
$\Gpd_RM\leq \dim R$.
\end{theorem}

\section{Gorenstein injective dimension}

Let $R$ be a local Noetherian ring an let ${\bf D}$ denote the
dualizing complex of $\hat{R}$. Let $B(\hat{R})$ denote the full
subcategory of ${\bf D}_b(\hat{R})$, consisting of those complexes
$X$ for which ${\bf R}\Hom_{\hat{R}}({\bf D},X)\in {\bf
D}_b(\hat{R})$ and the canonical morphism $$\tau_X:{\bf
D}\otimes_{\hat{R}}^{{\bf L}}{\bf R}\Hom_{\hat{R}}({\bf D},X)\lo
X,$$
is an isomorphism, see [{\bf 4}, 3.1].\\
\\
Now, we define $B'(R)$ to be the class of all $R$-modules $M$ such
that $\Hom_R(\hat{R},M)\in B(\hat{R})$.

In the Theorem 4.8, we want to characterize Gorenstein injective
modules in terms of the class $B'(R)$. To prove Theorem 4.8, we
need the following results.

\begin{definition} (See [{\bf 6}, Definition 5.10])
For every $R$-module $M$, we show  the large restricted injective
dimension by $\Ed_RM$ and define
$$\Ed_RM=\sup\{i\in \mathbb{N}_0\mid
\exists L\in \overline{F(R)}\mid \Ext_R^ i(L,M)\neq 0\}.$$
\end{definition}

\begin{theorem} ({\bf Dimension inequality}) Let $R$ be a Noetherian
ring of finite Krull dimension. For every $R$-module $M$, we have
the following inequality: $$\Ed_RM\leq \Gid_RM\leq \id_RM.$$
\end{theorem}

{\bf Proof.} Every injective module is Gorenstein injective, and so
$\Gid_RM\leq \id_RM$. We can assume that $\Gid_RM$ is finite. We
show that $\Ed_RM\leq \Gid_RM$ by induction on $\Gid_RM=n$. First
assume that $n=0$. By [{\bf 5}, Lemma 2.2], for every $R$-module $L
\in \overline{P(R)}$ and $i>0$, we have $\Ext_R^i(L,M)=0$. This
means that $\Ed_RM\leq 0$, and so that the result holds. Now, let
$n>0$. Using the Gorenstein injective version of [{\bf 13},
Proposition 2.18], there exists exact sequence $0\lo M\lo T\lo K\lo
0$ such that $T$ is injective $R$-module and $\Gid_RK=n-1$. By
induction, we have $\Ed_RK\leq \Gid_RK=n-1$, and so
$\Ext_R^j(L,K)=0$ for all $L\in \overline{F(R)}$ and all $j> n-1$.
For each $i>n$ and each $L\in \overline{F(R)}$, we have the
following exact sequence
$$0=\Ext_R^{i-1}(L,K)\lo \Ext_R^i(L,M)\lo \Ext_R^i(L,T)=0.$$
So $\Ed_RM\leq n=\Gid_RM$. This ends the proof.  $\Box$

By Theorem 4.2, every Gorenstein injective $R$-module over a
Noetherian ring of finite Krull dimension is strongly cotorsion
(see [{\bf 16}, Definition 5.4.1]). The following example shows
that there exists an $R$-module with finite Gorenstein injective
dimension over a regular local ring which is not cotorsion.

\begin{example} Let $R$ be a regular local ring of Krull
dimension one which is not complete. By [{\bf 1}, Lemma 3.3],
$\Hom_R(\hat{R},R)=0$. So, $\hat{R}$ is not a projective
$R$-module. Therefore, $\pd_R(\hat{R})=1$ and consequently there
exists an $R$-module $M$ such that $\Ext_R^ 1(\hat{R},M)\neq 0$.
On the other hand, $\id_RM\leq 1$. So, $M$ is an $R$-module with
finite Gorenstien injective dimension which is not cotorsion.
\end{example}

\begin{proposition} Let $R$ be a local Noetherian ring and $M$ an $R$-module.\\
(i) If $M$ is a Gorenstein injective $R$-module, then
$\Hom_R(\hat{R},M)$ is Gorenstein injective as an $\hat{R}$-module.\\
(ii) If $M$ is a Gorenstein injective $R$-module, then $M\in
B'(R)$.
\end{proposition}

{\bf Proof.} (i) Let
$${\bf X}=\ldots\lo E_2\lo E_1\lo E_0\overset{\rho^0}\lo G^0\lo G^1\lo
\ldots$$ be an exact sequence of injective $R$-modules such that
$\Hom_R(I,{\bf X})$ is exact for every injective $R$-modules $I$
with $\ker \rho^0=M$. If $$0\lo G{''}\lo E\lo G'\lo 0$$ is an
exact sequence such that $G'$, $G{''}$ are Gorenstein injective
and $E$ is injective, then Theorem 4.2 yields the short exact
sequence, $$0\lo \Hom_R(\hat{R},G{''})\lo \Hom_R(\hat{R},E)\lo
\Hom_R(\hat{R},G')\lo 0.$$ Hence, we obtain the exact sequence
$${\bf Y}=\ldots\lo \Hom_R(\hat{R},E_1)\lo
\Hom_R(\hat{R},E_0)\overset{\Hom_R(\hat{R},\rho^0)}\lo
\Hom_R(\hat{R},G^0)\lo \ldots $$ of $\hat{R}$-modules and
$\hat{R}$-homomorphisms in which $\ker(
\Hom_R(\hat{R},\rho^0))\cong \Hom_R(\hat{R},M)$. On the other
hand, if $E$ is an injective $R$-module, we can conclude that
$\Hom_R(\hat{R},E)$ is injective as an $\hat{R}$-module, because
$\Hom_{\hat{R}}(-,\Hom_R(\hat{R},E))\cong
\Hom_R(-\otimes_{\hat{R}}\hat{R},E)$. It is enough to show that
$\Hom_{\hat{R}}(E',{\bf Y})$ is exact, for all injective
$\hat{R}$-modules $E'$. This follows from the following
isomorphisms of complexes $$\Hom_{\hat{R}}(E',{\bf Y})\cong
\Hom_{\hat{R}}(E',\Hom_R(\hat{R},{\bf X}))\cong \Hom_R(E',{\bf
X})$$ and the fact that every injective $\hat{R}$-module is also
injective as an $R$-module.

(ii) Let $M$ be Gorenstein injective. By (i), $\Hom_R(\hat{R},M)$
is Gorenstein injective $\hat{R}$-module. Hence, by [{\bf 5},
Theorem 4.4], $\Hom_R(\hat{R},M)\in B({\hat{R}})$, and so $M\in
B'(R)$, by the definition. $\Box$

\begin{proposition} An $R$-module $M$ is Gorenstein injective if
and only if $\Ext_R^ i(E,M)=0$ for all injective $R$-modules $E$
and for all $i>0$ and there exists an exact sequence
$${\bf X}=\ldots\lo E_2\lo E_1\lo E_0\lo M\lo 0$$ of $R$-modules
and $R$-homomorphisms with $E_i$ is injective $R$-module for all
$i\geq 0$, such that $\Hom_R(E,{\bf X})$ is exact for all
injective $R$-modules $E$ (i.e. $M$ admits a proper left injective
resolution).
\end{proposition}

{\bf Proof.} It is the dual version of [{\bf 13}, Proposition 2.3]
and we leave the proof to the reader. $\Box$

\begin{lemma} (i) Let $R$ be a local Noetherian ring
and $M$ a cotorsion $R$-module such that $M$ belongs to $B'(R)$.
Then there exists an epimorphism $L\lo M$ with
$\id_R(L)<\infty$.\\
(ii) Let $R$ be a Noetherian ring and $\varphi:L\lo M$ an
$R$-epimorphism with $\id_R(L)<\infty$ and $\Ext_R^ i(N,M)=0$ for
all injective $R$-modules $N$ and all $i>0$. Then there exists an
epic $\overline{I(R)}$-precover $E\lo M$, in which $E$ is
injective.
\end{lemma}

{\bf Proof.} (i) Since $M$ belongs to $B'(R)$, then
$\Hom_R(\hat{R},M))$ belongs to $B(\hat{R})$. So,
$\Hom_R(\hat{R},M))$ has finite Gorenstein injective dimension as
an $\hat{R}$-module by [{\bf 5}, Theorem 4.4]. By [{\bf 5}, Lemma
2.18], There are an $\hat{R}$-module $L$ and an
$\hat{R}$-epimorphism $L\lo \Hom_R({\hat{R}},M)$ such that
injective dimension  of $L$ as an $\hat{R}$-module is finite.
Since every injective $\hat{R}$-module is injective as an
$R$-module, injective dimension of $L$ as an $R$-module is finite.
Consider the following exact sequence
$$0\lo R\lo \hat{R}\lo \hat{R}/R\lo 0,$$ that yields the following exact sequence
$$\Hom_R(\hat{R},M))\lo \Hom_R(R,M))\lo \Ext_R^ 1(\hat{R}/R,M).$$
On the other hand, since $\hat{R}/R$ is a flat $R$-module and $M$
is a cotorsion $R$-module, $\Ext_R^ 1(\hat{R}/R,M)=0$. So, the
natural $R$-homomorphism $\Hom_R(\hat{R},M)\lo M$ is epic. The
result follows.

(ii) By [{\bf 16}, Theorem 2.4.3], there exists an $I(R)$-precover
$f:E\lo M$. We claim that $f$ is an $\overline{I(R)}$-precover.
Let $\varphi':L'\lo M$ be an $R$-homomorphism such that
$\id_R(L')<\infty$. Consider an exact sequence
$$0\lo L'\overset{g}\lo E'\lo K\lo 0$$ such that $E'$ is an
injective $R$-module. It is clear that injective dimension of $K$
is finite. By induction on injective dimension, we can deduce from
assumption that $\Ext_R^ 1(K,M)$ is zero. We obtain the following
exact sequence
$$0\lo \Hom_R(K,M)\lo \Hom_R(E',M)\lo \Hom_R(L',M)\lo \Ext_R^
1(K,M)=0.$$  Hence, we conclude that there exists an
$R$-homomorphism $\psi:E'\lo M$ such that $\varphi'=\psi g$. On
the other hand, since $f$ is an $I(R)$-precover, there exists an
$R$-homomorphism $h:E'\lo E$ such that $\psi=f h$. Hence, there
exists an $R$-homomorphism $h g:L'\lo E$ such that $f (h
g)=\varphi'$. It therefore follows that $f$ is an
$\overline{I(R)}$-precover. Consequently $f$ is epic, because
$\varphi$ is epic. $\Box$

\begin{lemma} Let $(R,\fm)$ be a local Noetherian ring, $M$ a cotorsion
$R$-module, and $K$ a cotorsion $\hat{R}$-module. Then\\
(i) $\Ext_R^ i(F,M)=0$ for all flat $R$-modules $F$ and all $i>
0$.\\
(ii) $K$ is cotorsion as an $R$-module.\\
(iii) For all $j> 0$, $\Ext_R^ j(E,M)=0$ for all injective
$R$-modules $E$ if and only if $\Ext_{\hat{R}}^
j(I,\Hom_R(\hat{R},M))=0$ for all injective $\hat{R}$-modules $I$.
\end{lemma}

{\bf Proof.} (i) See the proof of [{\bf 16}, Proposition 3.1.2].\\
(ii) Suppose $F$ is a flat $R$-module and ${\bf P_{\bullet}}\lo F$
a projective resolution of $F$. For all $i>0$, we have
\begin{equation*}
\begin{split}
\Ext_R^i(F,K)& \cong H^ i(\Hom_R({\bf
P_{\bullet}},K))\\ &\cong H^ i(\Hom_{\hat{R}}({\bf P_{\bullet}}\otimes_R\hat{R},K))\\
&\cong \Ext_{\hat{R}}^ i(F\otimes_R\hat{R},K).
\end{split}
\end{equation*}
The last isomorphism comes from the fact that $K$ is a cotorsion
$\hat{R}$-module and $F\otimes_R\hat{R}$ is flat as an
$\hat{R}$-module for all flat
$R$-modules $F$. This ends the proof of (ii).\\
(iii) Suppose $L$ is an $\hat{R}$-module and ${\bf F_{\bullet}}\lo
L$ is a free resolution of $L$, considered as an $\hat{R}$-module.
For every $j> 0$, we have
\begin{equation*}
\begin{split}
\Ext_{\hat{R}}^j(L,\Hom_R(\hat{R},M))& \cong H^
j(\Hom_{\hat{R}}({\bf F_{\bullet}},\Hom_R(\hat{R},M))\\&\cong H^
j(\Hom_R({\bf F_{\bullet}}\otimes_{\hat{R}}\hat{R},M)) \\&\cong H^
i(\Hom_R({\bf F_{\bullet}},M)) \\&\cong \Ext_R^ j(L,M)
\end{split}
\end{equation*}
The last isomorphism follows from  the fact that $M$ is cotorsion
and every flat $\hat{R}$-module is flat as an $R$-module.\\
$\Rightarrow)$ We know that every injective $\hat{R}$-module is
injective as an $R$-module. So, the result follows from the above
isomorphism.\\
$\Leftarrow)$ By assumption, it is easy to see that
$$\Ext_{\hat{R}}^ i(N,\Hom_R(\hat{R},M))=0,$$ for all
$\hat{R}$-modules $N$ of finite injective dimension and all $i>0$.
Let $E$ be an injective $R$-module and let $\Hom_R(-,E_R(R/\fm))$
denote by $(-)^{\vee}$. From the natural monomorphism $E\lo
(E^{\vee})^{\vee}$, we conclude that $E$ is a direct summand of
$(E^{\vee})^{\vee}$. So, it is enough to show that $\Ext_R^
i((E^{\vee})^{\vee},M)=0$ for all $i>0$. Since, by Lemma 2.7,
$\id_{\hat{R}}((E^{\vee})^{\vee})<\infty$, the result follows from
the above isomorphism. $\Box$

\begin{theorem} Let $R$ be a local Noetherian ring and
$M$ an $R$-module. Then the following conditions are equivalent:\\
(i) $M$ is Gorenstein injective.\\
(ii) $M$ is cotorsion and $\Hom_R(\hat{R},M)$ is Gorenstein
injective as an $\hat{R}$-module.\\
(iii) $M\in B'(R)$, $M$ is cotorsion and $\Ext_R^ i(E,M)=0$ for
all injective $R$-modules $E$ and all $i>0$.
\end{theorem}

{\bf Proof.} $(i) \Rightarrow (ii)$ This follows from Theorem 4.2
and Proposition 4.4.

$(ii)\Rightarrow (iii)$ By [{\bf 5}, Theorem 4.4],
$\Hom_R(\hat{R},M)$ belongs to $B(\hat{R})$, and so $M$ belongs to
$B'(R)$. Also, Proposition 4.5 implies that
$$\Ext_{\hat{R}}^ i((I,\Hom_R(\hat{R},M))=0$$ for all
injective $\hat{R}$-modules $I$ and all $i>0$. The result follows
from Lemma 4.7 (iii).

$(iii) \Rightarrow (i)$ In view of Proposition 4.5, it is enough
to show that $M$ admits a proper left injective resolution. It
follows from Lemma 4.6 (i) and (ii) that there exists an exact
sequence
$$0\lo B\lo E\overset{f} \lo M\lo 0$$ such that $f$ is
an $\overline{I(R)}$-precover and $E$ an injective $R$-module. It
is enough to show that $B$ satisfies the given assumptions on $M$.

Let $I$ be an injective $R$-module. It is easy to deduce from the
above exact sequence that $\Ext_R^ i(I,B)=0$ for all $i\geq 2$.
Also, we have the following exact sequence
$$\Hom_R(I,E)\lo \Hom_R(I,M)\lo \Ext_R^ 1(I,B)\lo \Ext_R^
1(I,E)=0.$$ On the other hand, $\Hom_R(I,f)$ is epimorphism. So
$\Ext_R^ 1(I,B)=0$.

Now, we prove that $B$ is a cotorsion $R$-module. In view of
assumption and Lemma 4.7, we conclude that
$$\Ext_{\hat{R}}^i(I,\Hom_R(\hat{R},M))=0$$ for all injective
$\hat{R}$-modules $I$ and all $i>0$. On the other hand, $M\in
B'(R)$ implies that $\Hom_R(\hat{R},M)\in B(\hat{R})$. Therefore,
by [{{\bf 5}, Lemma 4.7}], $\Hom_R(\hat{R},M)$ is Gorenstein
injective as an $\hat{R}$-module. Hence, we have an exact sequence
$$0\lo K\lo E'\lo \Hom_R(\hat{R},M)\lo 0,$$ of $\hat{R}$-modules
and $\hat{R}$-homomorphism such that $E'$ is an injective and $K$
is a Gorenstein injective $\hat{R}$-module. By Theorem 4.2, $K$ is
a cotorsion $\hat{R}$-module. Lemma 4.7 implies that $K$ is
cotorsion as an $R$-module. Now, let $\varphi:\Hom_R(\hat{R},M)\lo
M$ be the natural $R$-homomorphism. Consider the following diagram

\begin{equation*} \setcounter{MaxMatrixCols}{11}
\begin{matrix}
&0\lo K&{\lo} E'\lo & \Hom_R(\hat{R},M)& \lo &0
\\
 & & & \mapdown{\varphi}
\\  &0\lo B& \lo E\overset{f}{\lo} & M& \lo &0.
\end{matrix}
\end{equation*}
Since $E'$ is an injective  $R$-module and $f:E\lo M$ is an
$\overline{I(R)}$-precover, there exists an $R$-homomorphism
$\psi:E'\lo E$ such that the following diagram is commutative.

\begin{equation*} \setcounter{MaxMatrixCols}{11}
\begin{matrix}
&0\lo K&{\lo} E'\lo & \Hom_R(\hat{R},M)& \lo &0
\\
 & & \mapdown{\psi} & \mapdown{\varphi}
\\  &0\lo B& \lo E\overset{f}{\lo} & M& \lo &0.
\end{matrix}
\end{equation*}
It is easy to see that there exists an $R$-homomorphism
$\theta:K\lo B$ such that the following diagram is commutative.

\begin{equation*} \setcounter{MaxMatrixCols}{11}
\begin{matrix}
&0\lo &K&{\lo} E'\lo & \Hom_R(\hat{R},M)& \lo &0
\\
 & & \mapdown{\theta} & \mapdown{\psi} & \mapdown{\varphi}
\\  &0\lo &B& \lo E\overset{f}{\lo} & M& \lo &0
\end{matrix}
\end{equation*}
Suppose $F$ is a flat $R$-module. Then we obtain the following
commutative diagram
\begin{equation*} \setcounter{MaxMatrixCols}{11}
\begin{matrix}
& \Hom_R(F,\Hom_R(\hat{R},M))& \overset{\beta} {\lo} & \Ext_R^ 1
(F,K)  & \lo 0&
\\
& \mapdown{\Hom_R(F,\varphi)} &&\mapdown{\theta_1}&&(*)
\\& \Hom_R(F,M)& \overset{\delta} {\lo} & \Ext_R^ 1
(F,B)  & \lo 0.&
\end{matrix}
\end{equation*}
The natural exact sequence $$0\lo R\lo \hat{R}\lo \hat{R}/R\lo
0,$$ yields the exact sequence  $$0\lo \Hom_R(\hat{R}/R,M)\lo
\Hom_R(\hat{R},M)\overset{\varphi}\lo M\lo 0,$$ because $M$ is a
cotorsion $R$-module and $\hat{R}/R$ is a flat $R$-module. Thus,
we obtain the following exact sequence\\ $0\lo
\Hom_R(F,\Hom_R(\hat{R}/R,M))\lo
\Hom_R(F,\Hom_R(\hat{R},M))\overset{\Hom_R(F,\varphi)}\lo
\Hom_R(F,M)\lo $ $$\lo \Ext_R^ 1(F,\Hom_R(\hat{R}/R,M)).$$ Since
$M$ is a cotorsion and $\hat{R}/R$ is a flat $R$-module, $$\Ext_R^
1(F,\Hom_R(\hat{R}/R,M))\cong \Ext_R^ 1(F\otimes_R{\hat{R}/R},M))
.$$ On the other hand, $F\otimes_R{\hat{R}/R}$ is a flat
$R$-module, so $\Ext_R^ 1(F\otimes_R{\hat{R}/R},M))$ is zero
$R$-module. Therefore $\Hom_R(F,\varphi)$ is an epimorphism. By
(*), $\theta_1 \beta$ is epic and so $\theta_1$ is epic. Thus,
since $K$ is a cotorsion $R$-module, $\Ext_R^ 1(F,B)$ is the zero
module. This means that $B$ is cotorsion.

Now, we apply the functor $\Hom_R(\hat{R},-)$ on the following
exact sequence $$0\lo B\lo E\lo M\lo 0,$$ and obtain the exact
sequence $$0\lo \Hom_R(\hat{R},B)\lo \Hom_R(\hat{R},E)\lo
\Hom_R(\hat{R},M)\lo 0.$$ It is easy to see that
$\Hom_R(\hat{R},E)$ is an injective $\hat{R}$-module. Since
$\Hom_R(\hat{R},M)$ is Gorenstein injective as an
$\hat{R}$-module, by [{\bf 13}, theorem 2.25], $\Hom_R(\hat{R},B)$
has finite Gorenstein injective dimension. So, it follows from
[{\bf 5}, Theorem 4.4] that $B\in B'(R)$. This ends the proof. $\Box$ \\

The following example shows that the dual version of Theorem 3.4
is not true.

\begin{example} Let $R$ be a non-complete local Noetherian domain which is
not Gorenstein. By [{\bf 14}, Theorem 2.1], $\Gid_R(R)=\infty$. On
the other hand, by [{\bf 1}, Lemma 3.3], $\Hom_R(\hat{R},R)=0$. So
$R$ has infinite Gorenstein injective dimension as an $R$-module
but $R\in B'(R)$.
\end{example}

\begin{theorem} Let $(R,\fm)$ be a local Noetherian ring of dimension $d$
and $\Ext_R^i(\hat{R},M)=0$ for all $i>0$. Then the Gorenstein
injective dimension of $M$ is finite if and only if $M$ belongs to
$B'(R)$. In particular, if $M\in B'(R)$ then $\Gid_R(M)\leq d$.
\end{theorem}
{\bf Proof.} $\Rightarrow )$ Let $\Gid_RM=t$ and $$0\lo M\lo
G^0\lo G^1\lo G^2\lo \ldots\lo G^t\lo 0$$ be an exact sequence
such that $G^i$ is Gorenstein injective for all $0\leq i\leq t$.
Using hypothesis, we obtain the following exact sequence
$$0\lo \Hom_R(\hat{R},M)\lo \Hom_R(\hat{R},G^0)\lo
\ldots\lo \Hom_R(\hat{R},G^t)\lo 0.$$ By Proposition 4.4 (i),
$\Gid_{\hat{R}}(\Hom_R(\hat{R},M))$ is finite as an
$\hat{R}$-module and so by [{\bf 5}, Theorem 4.4],
$\Hom_R(\hat{R},M)$ belongs to $B(\hat{R})$. The assertion
follows from the definition.

$\Leftarrow )$ Since $M$ belongs to $B'(R)$, $\Hom_R(\hat{R},M)$
belongs to $B(\hat{R})$. Now, by using [{\bf 5}, Theorem 4.4], the
Gorenstein injective dimension of $\Hom_R(\hat{R},M)$ is finite as
an $\hat{R}$-module. By [{\bf 13}, Theorem 2.29],
$\Gid_{\hat{R}}(\Hom_R(\hat{R},M))\leq \FID(R)$, where
$\FID(R)=\sup \{\ \id_R(M)| M \ is \ an \ R-module \ of \ finite \
injective \ dimension \}$. It is known that
$\id_R(N)=\fd_R(\Hom_R(N,E_R(R/{\fm})))$, for all $R$-modules $N$.
So, we have $\Gid_{\hat{R}}(\Hom_R(\hat{R},M))\leq d$.

Consider the following exact sequence $$0\lo M\lo E^ 0\lo E^ 1\lo
\ldots \lo E^{d-1}\lo L\lo 0,$$ of $R$-modules and
$R$-homomorphisms such that $E^ i$ is injective $R$-module for all
$0\leq i\leq d-1$. We have the following exact sequence,
$$0\lo \Hom_R(\hat{R},M)\lo \ldots \lo \Hom_R(\hat{R},E^{d-1})
\lo \Hom_R(\hat{R},L)\lo 0.$$ So, by [{\bf 13}, Theorem 2.22],
$\Hom_R(\hat{R},L)$ is a Gorenstein injective $\hat{R}$-module. On
the other hand, for any flat $R$-module $F$ and any $i> 0$, we
have
$$\Ext_R^ i(F,L)\cong \Ext_R^{i+d}(F,M).$$ Therefore, $\Ext_R^ i(F,L)$
is zero for all $i> 0$, because the projective dimension of $F$ is
less than $d+1$. So, $L$ is cotorsion. It therefore follows from
Theorem 4.8 that $L$ is a Gorenstein injective $R$-module. Thus,
$\Gid_R(M)\leq d$. $\Box$ \\

{\bf Acknowledgement.} In the first version of the manuscript, we
proved our main results for local Cohen-Macaulay rings. We would
like to thank Lars Winther Christensen for pointing out that our
main results can be extended to non Cohen-Macaulay case.


\end{document}